\documentstyle{amltd2004}
\input boxedeps.tex 
\SetepsfEPSFSpecial 
\HideDisplacementBoxes
\def\figin#1#2{
$$
 {\BoxedEPSF{#1.eps scaled
#2}%
}%
$$
\medbreak\noindent}
\begin{document}
\annalsline{155}{2002}
\received{May 19, 2001}
\startingpage{295}
\def\bye{\end{document}}
 \font\tenrm=cmr10

\def\joinrel{\mathrel{\mkern-4mu}}
\def\relbar{\mathrel{\smash-}}
\def\lrar{\relbar\joinrel\relbar\joinrel\relbar\joinrel\rightarrow}
\def\ritem#1{\item[{\rm #1}]}

\input amssym.def
\input amssym.tex

\catcode`\@=11
\font\twelvemsb=msbm10 scaled 1100
\font\tenmsb=msbm10
\font\ninemsb=msbm10 scaled 800
\newfam\msbfam
\textfont\msbfam=\twelvemsb  \scriptfont\msbfam=\ninemsb
  \scriptscriptfont\msbfam=\ninemsb
\def\msb@{\hexnumber@\msbfam}
\def\Bbb{\relax\ifmmode\let\next\Bbb@\else
 \def\next{\errmessage{Use \string\Bbb\space only in math
mode}}\fi\next}
\def\Bbb@#1{{\Bbb@@{#1}}}
\def\Bbb@@#1{\fam\msbfam#1}
\catcode`\@=12

 \catcode`\@=11
\font\twelveeuf=eufm10 scaled 1100
\font\teneuf=eufm10
\font\nineeuf=eufm7 scaled 1100
\newfam\euffam
\textfont\euffam=\twelveeuf  \scriptfont\euffam=\teneuf
  \scriptscriptfont\euffam=\nineeuf
\def\euf@{\hexnumber@\euffam}
\def\frak{\relax\ifmmode\let\next\frak@\else
 \def\next{\errmessage{Use \string\frak\space only in math
mode}}\fi\next}
\def\frak@#1{{\frak@@{#1}}}
\def\frak@@#1{\fam\euffam#1}
\catcode`\@=12

\newcommand{\Curve}{{\cal C}}
\renewcommand{\O}{{\cal O}}
\newcommand{\FF}{{\cal F}}
\newcommand{\EE}{{\cal E}}
\newcommand{\F}{{\Bbb F}}
\newcommand{\Fp}{{{\Bbb F}_p}}
\newcommand{\Fq}{{{\Bbb F}_q}}
\newcommand{\Z}{{\Bbb Z}}
\newcommand{\Q}{{\Bbb Q}}
\newcommand{\C}{{\Bbb C}}
\renewcommand{\P}{{\Bbb P}}
\newcommand{\A}{{\Bbb A}}
\newcommand{\p}{{\frak p}}
\newcommand{\n}{{\frak n}}
\newcommand{\into}{\hookrightarrow}
\newcommand{\isoto}{\tilde{\to}}
\newcommand{\tensor}{\otimes}
\newcommand{\compose}{\circ}
\def\nodiv{\mathrel{\mathchoice{\not|}{\not|}
{\kern-.2em\not\kern.2em|}{\kern-.2em\not\kern.2em|}}}
\def\tr{\mathop{\rm Tr}}
\def\N{\mathop{\rm N}}
\def\ord{\mathop{\rm ord}}
\def\rk{\mathop{\rm Rank}}
\def\aut{\mathop{\rm Aut}}
\def\gal{\mathop{\rm Gal}}
\def\spec{\mathop{\rm Spec}}

\newcommand{\Fbar}{{\overline{{\Bbb F}}}}
\newcommand{\Fpbar}{{\overline{{\Bbb F}}_p}}
\newcommand{\Fqbar}{{\overline{{\Bbb F}}_q}}
\newcommand{\Fpf}{{{\Bbb F}}_{p^f}}
\newcommand{\Qbar}{{\overline{{\Bbb Q}}}}
\newcommand{\EEbar}{{\overline{\EE}}}
\newcommand{\kbar}{{\overline{k}}}
\newcommand{\ratto}{{\dashrightarrow}}

\title{Elliptic curves with large
rank \\ over function fields} 
\shorttitle{Elliptic curves with large rank} 

  \acknowledgements{This paper is based upon work supported by the National
Science Foundation under Grant No. DMS-0070839.}
\author{Douglas Ulmer}
\institutions{University of Arizona, Tucson,
AZ\\
{\eightpoint {\it E-mail address\/}: ulmer@math.arizona.edu}}

\bigbreak\centerline{\bf Abstract}
\medbreak

We produce explicit elliptic curves over $\Fp(t)$ whose Mordell-Weil
groups have arbitrarily large rank.  Our method is to prove
the conjecture of Birch and Swinnerton-Dyer for these curves (or
rather the Tate conjecture for related elliptic surfaces) and then
use zeta functions to determine the rank.  In contrast to earlier
examples of Shafarevitch and Tate, our curves are not isotrivial.

Asymptotically these curves have maximal rank for their conductor.
Motivated by this fact, we make a conjecture about the growth of ranks
of elliptic curves over number fields.
 
\section{Introduction} 

 1.1. 
Let $K$ be a field and consider elliptic curves defined over $K$.  A
natural question is whether there exist elliptic curves over $K$ with
arbitrarily large Mordell-Weil rank.  In other words, for every $r$ is
there an $E$ such that $E(K)$ contains at least $r$ independent points
of infinite order?  The general expectation seems to be that such
curves exist for any field $K$ which is not algebraic over a finite
field.
 
\demo{{\rm 1.2}} 
For fields of characteristic zero, it obviously suffices to treat the
case $K=\Q$.  Here the question is open and it seems to be quite
difficult to produce examples with large rank.  At this writing, the
largest known rank is 24 (\cite{MM00}) and the largest
proven analytic rank (i.e., order of vanishing of $L$-series at $s=1$)
is 3 (\cite{GZ86}).
\enddemo

1.3. 
For fields of characteristic $p$, it suffices to consider the rational
function field $K=\F_p(t)$.  In~\cite{TS67},
Shafarevitch and Tate produced elliptic curves over $K$ of arbitrarily
large rank.  They considered a supersingular curve $E_0$ defined over
$\F_p$ (viewed as a curve $E$ over $K$ in the obvious way) and showed
that that there are quadratic extensions $L/K$ such that the Jacobian
of the curve over $\F_p$ attached to $L$ has a large number of factors
isogenous to $E_0$ over~$\F_p$.  This implies that the quadratic twist
of $E$ by $L$ has large rank.

\demo{{\rm 1.4}} 
The examples of Shafarevitch and Tate are ``isotrivial,'' i.e., after
a finite extension ($L/K$ in fact), they become isomorphic to elliptic
curves defined over $\F_p$.  (Equivalently, their $j$-invariants lie
in $\F_p$.)  There is no analog of this property for an elliptic curve
over $\Q$ and so it is not clear whether their examples provide
evidence for the question over $\Q$.  (On the other hand,
isotriviality makes sense for fields like $\Q(t)$, and it is
conceivable that the arguments of Tate and Shafarevitch might be
generalized to this context.) 
\enddemo 
The aim of this paper is to produce elliptic curves over $K=\F_p(t)$
which are nonisotrivial ($j\not\in\F_p$) and which have arbitrarily
large rank.

\specialnumber{1.5} \proclaim{Theorem}\label{thm:intro}
Let $p$ be an arbitrary prime number{\rm ,} $\Fp$ the field of $p$ elements{\rm ,}
and $\Fp(t)$ the rational function field in one variable over $\Fp${\rm .} 
Let $E$ be the elliptic curve defined over $K=\F_p(t)$ by the
Weierstrass equation
$$
y^2+xy=x^3-t^d
$$
where $d=p^n+1$ and $n$ is a positive integer{\rm .}  Then
$j(E)\not\in\F_p${\rm ,} the conjecture of Birch and Swinnerton\/{\rm -}\/Dyer holds
for $E$ over $K${\rm ,} and the rank of $E(K)$ is at least $(p^n-1)/2n${\rm .}
\endproclaim

\demo{{\rm 1.6.} Remarks}
\begin{itemize}
\item[1.]
We give a simple expression for the exact rank of $E(K)$ in
Theorem~\ref{thm:main}. 
The $j$-invariant of $E$ is $t^{-d}(1-2^43^3t^d)^{-1}$.
\item[2.]
In fact we prove the conjecture of Birch and Swinnerton-Dyer for
$E$ over $\F_q(t)$ for $q$ any power of $p$, and we show that
$\rk E(\F_{p^{2n}}(t))=\rk E(\Fbar_p(t))=p^n$ if $6\nodiv d$ and
$p^n-2$ if $6|d$.
\item[3.]
In Section~\ref{s:Bounds} we explain that the curves in
Theorem~\ref{thm:intro} asymptotically have maximal ranks for their
conductor and we make a conjecture about ranks of elliptic curves
over number fields.
\item[4.]
The displayed Weierstrass equation also defines an elliptic curve over
$\Q(t)$.  It turns out that this curve has rank which is bounded
independently of $d$, even over $\Qbar(t)$.
\end{itemize}

\enddemo

1.7.
The proof of the theorem involves an appealing mix of geometry and
arithmetic.  We begin with the geometry: First, we construct an
elliptic surface $\EE\to\P^1$ over $\F_p$ whose generic fiber is
$E/K$.  The rank of the Mordell-Weil group $E(K)$ is closely related
to the rank of the N\'eron-Severi group of $\EE$, i.e., to curves on
$\EE$ up to algebraic equivalence.  Next, we define a dominant
rational map from a Fermat surface $F_d$ to $\EE$, which induces a
birational isomorphism between $\EE$ and a certain quotient
$F_d/\Gamma$ of the Fermat surface.  Thirdly, we carry out a fairly
detailed analysis of the geometry of this birational map.

Then comes the arithmetic: The Tate conjecture (on cycles and poles of
zeta functions) is known for Fermat surfaces and this allows us to
deduce it for $\EE$.  (The conjecture of Birch and Swinnerton-Dyer for
$E$ is equivalent to the Tate conjecture for $\EE$.)  Also, the
detailed analysis of the birational isomorphism between $F_d/\Gamma$
and $\EE$ allows us to express the zeta function of $\EE$ in terms of
that of $F_d$ which was calculated by Weil in terms of Gauss sums.
Finally, an explicit calculation of Gauss sums allows us to show that
the zeta function of $\EE$ has a large order pole at $s=1$, and
therefore $E(K)$ has large rank.

\demo{{\rm 1.8}} 
One of the key ideas of the proof, namely relating $\EE$ to a Fermat
surface, is due to Shioda.  In~\cite[\S 5]{Shi86}, he
exhibited (nonisotrivial) elliptic curves $E_n$ ($n\ge1$ an integer)
defined over $\F_p(t)$ (with $p\equiv3\; ({\rm mod}\, 4)$) which are related in a
similar way to Fermat surfaces.  This allowed him to calculate the
rank of $E_n$ over $\Fbar_p(t)$ and thus to show that it tends to
infinity with $n$.  (In fact, it is not difficult to see that Shioda's
curves achieve their full rank over $\F_{p^{2n}}(t)$; i.e.,
$E_n(\F_{p^{2n}}(t))=E_n(\Fbar_p(t))$.)  The methods of this paper can
be used to show that the rank of $E_n(\F_p(t))$ also tends to infinity
(although in general, $E_n(\F_p(t))$ has smaller rank than
$E_n(\Fbar_p(t))$).
\enddemo

1.9.
Our theorem says that $E(K)$ has large rank, but the proof does not
explicitly produce any points.  Exhibiting explicit points, and
computing invariants such as the height pairing, looks like an
interesting project.  Do these Mordell-Weil lattices have high
densities or other special properties?

\demo{{\rm 1.10}} 
It is a pleasure to thank Felipe Voloch for bringing Shioda's paper
\cite{Shi86} to my attention, Pavlos Tzermias for his help with
$p$-adic Gamma functions and the Gross-Koblitz formula, and Dinesh
Thakur for a number of useful remarks.
 
\vglue-8pt
\section{Invariants of $E$}\label{s:E}
\vglue-4pt

2.1.
We work in somewhat greater generality than
in the introduction.  Let $k$ be a perfect field (of characteristic
$p=0$ or a prime) and set $K=k(t)$.  Fix a positive integer $d$ which
is not divisible by $p$. Let $E_d$ be the elliptic curve over $K$ with
plane cubic model
\begin{equation}\label{eq:E}
y^2+xy=x^3-t^d. \speqnu{2.1.1}
\end{equation}
We will usually drop $d$ from the notation and write $E$ for $E_d$. 

\demo{{\rm 2.2}} 
Straightforward calculation shows that the discriminant of this model is
$\Delta=t^d(1-2^43^3t^d)$ and $j(E)=1/\Delta$.  Thus $E$ has good
reduction at all places of $K$ except $t=0$, the divisors of
$(1-2^43^3t^d)$ and possibly $t=\infty$.  

Applying Tate's algorithm (\cite{Tat75}), we see that at $t=0$, $E$
has split multiplicative reduction of type $I_d$ and all geometric
components of the special fiber are rational over $k$.  At places $v$
dividing $(1-2^43^3t^d)$, $E$ has multiplicative reduction of type
$I_1$; the tangent directions at the node are rational over
$k_v(\mu_4)$ where $k_v$ is the residue field at $v$ and $\mu_4$
denotes the 4th roots of unity.
\enddemo

\vglue-4pt
2.3.
The reduction type of $E$ at $t=\infty$ depends on $d\, ({\rm mod}\, 6)$ and on~$k$.  Write $d=6a-b$ where $0\le b<6$. 
Changing coordinates via $x=t^{2a}x'$, $y=t^{3a}y'$, and $t={t'}^{-1}$ we have the
model
\smallbreak
\centerline{$
{y'}^2+{t'}^ax'y'={x'}^3-{t'}^b.
$}
\smallbreak

Applying Tate's algorithm, we find the data in the following table.
$$
\begin{array}{|l||c|c|c|c|c|c|}
\hline
b&0&1&2&3&4&5\\
\hline
\mathrm{Reduction}&I&II&IV&I_0^*&IV^*&II^*\\
\hline
c&1&1&\begin{array}{l}
\strut
1\hbox{ if }\mu_4\not\subset k\\3\hbox{ if }\mu_4\subset k\end{array}&
\begin{array}{l}
\strut
2\hbox{ if }\mu_3\not\subset k\\4\hbox{ if }\mu_3\subset k\end{array}&
\begin{array}{l}
\strut
1\hbox{ if }\mu_4\not\subset k\\3\hbox{ if }\mu_4\subset k\end{array}&
1\\
\hline
n&1&1&\begin{array}{l}
\strut
2\hbox{ if }\mu_4\not\subset k\\3\hbox{ if }\mu_4\subset k\end{array}&
\begin{array}{l}
\strut
4\hbox{ if }\mu_3\not\subset k\\5\hbox{ if }\mu_3\subset k\end{array}&
\begin{array}{l}
\strut
5\hbox{ if }\mu_4\not\subset k\\7\hbox{ if }\mu_4\subset k\end{array}&
9\\
\hline
\end{array}
$$
Here $c$ is the number of geometric components of multiplicity 1 in
the special fiber which are $k$-rational and $n$ is the number of
irreducible components of the special fiber as a scheme over $k$.

The exponent of the conductor at $t=\infty$ is 0 if
$b=0$; 2 if $b>0$ and $p\nodiv6$; and $d+2$ if $p|6$.

\vglue-12pt
\section{Construction of $\EE$}\label{s:EE}
\vglue-6pt

3.1.
As before, we let $k$ be a perfect field of characteristic $p$
(possibly 0), $K=k(t)$, and $E$ the elliptic curve over $K$ defined by
Equation~\ref{eq:E}.  Our purpose in this section is to construct the
unique elliptic surface $\pi:\EE\to\P^1$ over $k$ such that $\EE$ is
regular and $\pi$ is proper, flat, and relatively minimal, with
generic fiber $E\to\spec K$.  We also relate the N\'eron-Severi group
of $\EE$ to the Mordell-Weil group $E(K)$.

\vglue-4pt
\demo{{\rm 3.2}} 
To that end, let $U$ be the closed subset of $\P^2\times\A^1$ over $k$
defined by the equation
 
\centerline{$
y^2z+xyz-x^3+z^3t^d=0
$}
\smallbreak\noindent 
where $x$, $y$, and $z$ are the coordinates on $\P^2$ and $t$ is the
coordinate on $\A^1$.  Similarly, define $U'\subset\P^2\times\A^1$ by
the equation
\smallbreak
\centerline{$
{y'}^2z'+{t'}^ax'y'z'-{x'}^3+{z'}^3{t'}^b=0
$}
\smallbreak\noindent where $d=6a-b$ with $0\le b<6$.  
\pagebreak

Let $W$ be the result of glueing $U\setminus\{t=0\}$ and
$U'\setminus\{t'=0\}$ via the identification
$([x',y',z'],t')=([t^{-2a}x,t^{-3a}y,z],t^{-1})$.  ($W$ stands for
``Weierstrass model.'')  Projection onto the $t$ or $t'$ coordinate
gives a proper, flat, and relatively minimal morphism $W\to\P^1$ whose
fibers are the naive reductions of $E$ at places of $K$.  However, $W$
may not be regular; there are singularities at the point $([0,0,1],0)$
in $U$ if $d>1$, and at $([0,0,1],0)$ in $U'$ if $b>1$.  We note that
$W$ is a local complete intersection and nonsingular in codimension
1, so is normal.
\enddemo

To resolve the singularity in $U$ we must blow up $W$ $\lfloor
d/2\rfloor$ times.  The fiber over $t=0$ is then a chain of $d$
rational curves meeting transversally; i.e., it is the fiber in the
N\'eron model of $E$ at $t=0$.  Tate's algorithm gives the recipe for
resolving the singularity in $U'$. The fiber over $t=\infty$ was
recorded in the previous section.  We denote by $\EE$ the result of
this desingularization.  We have a commutative diagram
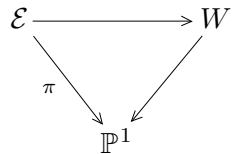
\figin{diag1}{1000}
\noindent 
where $\EE\to W$ is a birational isomorphism and the maps to $\P^1$ 
are proper, flat, and relatively minimal.  The map $\pi$ admits a
canonical section, the ``0-section,'' defined in the $U$ coordinates
by $t\mapsto ([0,1,0],t)$.  

\demo{{\rm 3.4}} 
Let $\EEbar=\EE\times_{\spec k}\spec\overline k$.  The N\'eron-Severi
group $NS(\EEbar)$ is by definition the group of divisors on $\EEbar$
modulo algebraic equivalence.  The ``theorem of the base'' asserts
that this is a finitely generated abelian group.  The N\'eron-Severi
group of $\EE$ is by definition the image of the group of divisors on
$\EE$ in $NS(\EEbar)$.  In other words, it is the group of ($k$-rational)
divisors on $\EE$ modulo algebraic equivalence over $\kbar$.  (If $k$
is finite, $NS(\EE)$ is the set of $\gal(\kbar/k)$-invariant elements
of $NS(\EEbar)$.)
\enddemo

3.5.
Let $L\subset NS(\EE)$ be the subgroup generated by the the class of
the 0-section and classes of divisors supported in fibers of $\pi$.  The
formula of Tate and Shioda (\cite{Tat66}, \cite{Shi72}) says that 
\begin{eqnarray*}
\rk E(K)&=&\rk NS(\EE)-\rk L\\
&=&\rk NS(\EE) - 2 -\sum_v(n_v-1)
\end{eqnarray*}
where the sum is over all closed points of $\P^1_k$
and $n_v$ is the number of irreducible components of the fiber at $v$.

Although they are not strictly necessary for our purposes, the
following remarks may help to clarify this formula.
We have an exact sequence
$$
0\to L\to NS(\EE)\to E(K)\to 0
$$
where the map $NS(\EE)\to E(K)$ can be defined as follows: given a divisor
on~$\EE$, take its intersection with the generic fiber $E$ and add the
resulting points in $E(K)$.  With some work, this can be shown to be
well-defined with kernel $L$.  

There is a section $s:E(K)\to NS(\EE)$ which sends a point in $E(K)$
to its closure in $\EE$.  But note that this section is not in general
a homomorphism.  (There is a canonical section after tensoring the
exact sequence with $\Q$ which makes $NS(\EE)\tensor\Q$ the orthogonal
direct sum, with respect to the intersection pairing, of $L\tensor\Q$
and $E(K)\tensor\Q$.)  The exact sequence is equivalent to the
assertion that that $s(E(K))$ is a set of coset representatives for
$L$ in $NS(\EE)$.  See \cite{MP86} for more details.

The Shioda-Tate formula will allow us to compute the rank of
$E(K)$ in terms of $NS(\EE)$, which we will eventually compute using
the Tate conjecture.

\demo{{\rm 3.6}} 
We will want to consider the situation for different degrees, so we
write $U_d$, $U'_d$, $W_d$ and $\EE_d$ for the surfaces considered in
this section; we then have similar definitions where $d$ is replaced
by 1.  Note that there is a finite morphism $W_d\to W_1$ defined in the
$U$ coordinates by $([x,y,z],t)\mapsto([x,y,z],t^d)$.
\enddemo

\section{Fermat surfaces}\label{s:Fd}

4.1.
In this section we will make a connection between Fermat surfaces and
the elliptic surface $\EE$.  

Let $F_d$ be the Fermat surface of degree $d$, i.e., the hypersurface
in $\P^3$ defined by
$$
x_0^d+x_1^d+x_2^d+x_3^d=0.
$$
We write $\mu_d$ for the group of $d^{\rm th}$ roots of unity in $\kbar$.
Let $G\subset\aut_\kbar(F_d)$ be the quotient of ${\mu_d}^4$ modulo a
diagonally embedded copy of $\mu_d$.  The action on $F_d$ is
$$
z\cdot x=
[\zeta_0,\zeta_1,\zeta_2,\zeta_3]\cdot
[x_0,x_1,x_2,x_3]=
[\zeta_0x_0,\zeta_1x_1,\zeta_2x_2,\zeta_3x_3].
$$
The canonical morphism $F_d\to F_1$
($[x_0,x_1,x_2,x_3]\mapsto[x_0^d,x_1^d,x_2^d,x_3^d]$) induces 
an isomorphism $F_d/G\cong F_1\cong\P^2$.

\demo{{\rm 4.2}} 
We define dominant rational maps $F_d\ratto W_d$ and $W_d\ratto F_1$
using the $U_d$ coordinates as follows:
$$
[x_0,x_1,x_2,x_3]\mapsto\left([(x_0x_1x_2)^d,(x_0^2x_2)^d,-x_1^{3d}],
\frac{x_0^3x_2^2x_3}{x_1^6}\right)
$$
and
$$
\left([x,y,z],t\right)\mapsto[y^2z,xyz,-x^3,z^3t^d].
$$
The rational map $W_d\ratto F_1$ factors as $W_d\to W_1\ratto F_1$ and
$W_1\ratto F_1$ is a birational isomorphism  with inverse
$$[x_0,x_1,x_2,x_3]\mapsto\left(\left[x_0x_1x_2,x_0^2x_2,-x_1^{3}\right],
x_0^3x_2^2x_3/x_1^6\right)$$ and so $W_d\ratto F_1$ has generic degree
$d$. Also, the composition $F_d\ratto W_d\ratto F_1$ is the canonical
morphism $F_d\to F_1$ which has degree $d^3$, so $F_d\ratto W_d$ has
generic degree $d^2$.
\enddemo

4.3.
Fix a primitive $d^{\rm th}$ root of unity $\zeta\in\kbar$ and let $\Gamma\subset G$ be the
subgroup generated by $[\zeta^2,\zeta,1,1]$ and
$[1,\zeta,\zeta^3,1]$, which can also be described as
$$
\Gamma=\left\{[\zeta_0,\zeta_1,\zeta_2,\zeta_3]\left|
\zeta_0^3\zeta_1^{-6}\zeta_2^2\zeta_3=1\right.\right\}.
$$
It is evident from the definitions that the rational map $F_d\ratto W_d$
factors through $F_d/\Gamma$.  Considering degrees, we see that the
induced map $F_d/\Gamma\ratto W_d$ is a birational isomorphism.  We
denote its inverse by $\varphi_d: W_d\ratto F_d/\Gamma$.

Note that since $F_d$ is regular, and thus normal, $F_d/\Gamma$ is
normal.

To summarize, we have the following commutative diagram, where the
vertical arrows are finite surjective morphisms, the horizontal arrows
are birational isomorphisms, and the diagonal arrows are dominant
rational maps.
\figin{diag2}{1000}
\noindent 
The surfaces $W_d$, $W_1$, and $F_d/\Gamma$ are normal, whereas $F_d$,
$\EE_d$, and $F_1$ are regular. 

\section{Analysis of $\varphi_d$}\label{s:geometry}

5.1.
We will eventually compute the zeta-function of $\EE_d$ in terms of that
of $F_d$.  In order to do this, we need some detailed geometric
information about $\varphi_d:W_d\ratto F_d/\Gamma$. Specifically, we
need to find explicit closed sets to be removed from $W_d$ and
$F_d/\Gamma$ so that $\varphi_d$ induces an isomorphism on the
remaining open sets.  Attacked directly, this computation could be
rather unpleasant, since we would need equations for $F_d/\Gamma$.

To avoid this problem, we prove a lemma which could be phrased
colloquially as saying that ``a rational map from a normal variety is
defined at a point if and only if its composition with a finite map
is.''

\specialnumber{5.2} \proclaim{Lemma}
Let $W${\rm ,} $X${\rm ,} $Y${\rm ,} and $Z$ be varieties over $k$ {\rm (}\/separated integral
schemes of finite type over $k${\rm )} and assume that $X$ is normal{\rm .}  Let
$g:X\ratto Y$ be a rational map and let $f:W\to X$ and $h:Y\to Z$ be
finite morphisms{\rm ,} with $f$ surjective{\rm .}
\vglue3pt
{\rm 1.} $g\compose f$ is defined at $w\in W$ if and only if $g$ is
defined at $f(w)${\rm .}
\vglue3pt {\rm 2.} $g$ is defined at $x\in X$ if and only if $h\compose g$ is
defined at $x${\rm .}
\endproclaim

{\it Proof}.
We may assume that all our varieties are affine, say
$W=\spec R$, $X=\spec S$, $Y=\spec T$, and $Z=\spec U$, and that we have
homomorphisms $g^*:T\to S[1/s]$ for some $s\in S$ and $f^*:S\to R$,
and $h^*:U\to T$.  The hypotheses imply that $R$, $S$, $T$, and $U$ are
domains, $S$ is integrally closed, $R$ is integral over $f^*(S)$, $T$
is integral over $h^*(U)$ and $f^*$ is injective.

To prove 1, we must show that $(g\compose f)^*$ factors through $R$ if
and only if $g^*$ factors through $S$.  The ``if'' direction is trivial.
For the converse, take $t\in T$ and write down an equation of
integrality for $f^*(g^*(t))$ over $f^*(S)$:
$$
f^*(g^*(t))^n+f^*(s_1)f^*(g^*(t))^{n-1}+\cdots+f^*(s_n)=0.
$$
Since $f^*$ is injective, this implies 
$$
g^*(t)^n+s_1g^*(t)^{n-1}+\cdots+s_n=0.
$$
But $S$ is integrally closed and $g^*(t)\in S[1/s]$, so $g^*(t)\in S$.

To prove 2, we must show that $(h\compose g)^*$ factors through $S$ if
and only if $g^*$ factors through $S$.  Again the ``if'' direction is
trivial.  For the converse, take $t\in T$ and write down an equation
of integrality over $h^*(U)$:
$$
t^n+h^*(u_1)t^{n-1}+\cdots+h^*(u_n)=0.
$$
Applying $g^*$ we have an equation of integrality for $g^*(t)$ over
$(h\compose g)^*(U)\subset S$.  Since $S$ is integrally closed and
$g^*(t)\in S[1/s]$, we have $g^*(t)\in S$.
\hfill\qed

\specialnumber{5.3} \proclaim{{C}orollary}
Consider a diagram of varieties
$$
\begin{array}{ccc}
 {\tilde X}&\hskip-.25in\stackrel{\tilde \phi}{- - \! \to}&\hskip-.25in {\tilde Y}\\
 {\scriptstyle\pi_X}\Big\downarrow\phantom{{\scriptstyle\pi_X}}&&\hskip-.28in
\phantom{\scriptstyle\pi_X}\Big\downarrow{\scriptstyle\pi_Y} \\
 X&\hskip-.25in\lower8pt\hbox{$\stackrel{\textstyle {- - \!\to}}{\scriptstyle\phi}$}&\hskip-.25in Y \end{array}
$$
where the vertical arrows are finite surjective morphisms{\rm ,} the
horizontal arrows are dominant rational maps and $\tilde X$ and $\tilde
Y$ are normal{\rm .}  If $V\subset X$ and $V'\subset Y$ are open subsets such
that $\phi$ induces a biregular isomorphism $\phi:V\isoto V'${\rm ,} then
$\tilde\phi$ induces a biregular isomorphism from $\tilde
V=\pi_X^{-1}(V)$ to $\tilde
{V'}=\pi_Y^{-1}(V')${\rm .}
\endproclaim  

\demo{Proof}
By the trivial half of part 1 of the lemma, $\phi\compose\pi_X$ is
defined on $\tilde V$, thus $\pi_Y\compose\tilde\phi$ is defined there
as well.  Part 2 of the lemma then implies that $\tilde\phi$ is
defined on $\tilde V$.  Similarly, $\tilde\phi^{-1}$ is defined on
$\tilde{V'}$.  Obviously $\tilde\phi(\tilde V)\subset \tilde{V'}$ and
$\tilde\phi^{-1}(\tilde{V'})\subset \tilde V$.  Since
$\tilde\phi^{-1}\compose\tilde\phi:\tilde V\to\tilde V$ and 
$\tilde\phi\compose\tilde\phi^{-1}:\tilde {V'}\to\tilde {V'}$
represent the rational maps $id_{\tilde X}$ and $id_{\tilde Y}$, they
must be the identity maps.\enddemo

\demo{{\rm 5.4}}
Now we apply the corollary to the diagram
$$
\begin{array}{ccc}
 {W_d}&\hskip-6pt\stackrel{ \varphi_d}{- - \! \to}&\hskip-6pt {F_d/\Gamma}\\
 \Big\downarrow &&\hskip-.28in
 \Big\downarrow \\ \noalign{\vskip5pt}
 W_1&\hskip-6pt\lower8pt\hbox{$\stackrel{\textstyle {- - \!\to}}{\scriptstyle\varphi_1}$}&\hskip-12pt F_1\, .
\end{array}
$$  
The rational map $\varphi_1$ is defined in the $U_1$ coordinates by
$$
\varphi_1([x,y,z],t)=[y^2z,xyz,-x^3,z^3t]
$$
and $\varphi_1^{-1}$ is defined by
$$
\varphi_1^{-1}([x_0,x_1,x_2,x_3])=
\left([x_0x_1x_2,x_0^2x_2,-x_1^3],\frac{x_0^3x_2^2x_3}{x_1^6}\right).
$$
We let $V\subset W_1$ be the subset of $U_1$ where $xyz\neq0$ and let
$V'\subset F_1$ be the subset where $x_0x_1x_2\neq0$.  It is easy to
see that $\varphi_1$ is defined on $V$, $\varphi_1^{-1}$ is defined on
$V'$, and they are inverse morphisms.   We conclude that $\varphi_d$
maps the subset of $U_d$ where $xyz\neq0$ isomorphically onto a certain
subset of $F_d/\Gamma$ which will be described in the next
subsection. 
\enddemo

5.5.
Consider the projection morphisms $F_d\to F_d/\Gamma\to F_1$.
The subset $x_0x_1x_2=0$ of $F_1$ has three irreducible components,
all lines.  Its inverse image in $F_d$ consists of three irreducible
curves (Fermat curves of degree $d$).  It follows that the inverse
image of $x_0x_1x_2=0$ in $F_d/\Gamma$ is also the union of three
irreducible curves.  Thus the open set
$V'\subset F_d/\Gamma$ of the preceding subsection is the complement
of the union of three irreducible curves.

\demo{{\rm 5.6}}
Now we consider the open subset $V\subset W_d$.  It is obtained from
$W_d$ by removing the subset where $t'=0$ from $U'_d$, as well as the
subset where $xyz=0$ from $U_d$.

The subset of $U'_d$ where $t'=0$ is an irreducible curve.

The subset $\{x=0\}=\{x=0=y^2z+z^3t^d\}$ of $U_d$ consists of the zero
section ($z=0$) and one or two other components, 1 if $d$ is odd or
$\mu_4\not\subset k$ and two if $d$ is even and $\mu_4\subset k$.  

The subset $\{y=0\}=\{y=0=x^3-z^3t^d\}$ of $U_d$ is irreducible if
$3\nodiv d$, it has two components if $3|d$ and $\mu_3\not\subset k$,
and it has three components if $3|d$ and $\mu_3\subset k$.

The subset where $z=0$ is contained in the subset where $x=0$.

In summary, $V$ is obtained from $W_d$ by removing a closed subset
which is a union of curves.  The number of irreducible components of
this union is
$$
1+\left\{
\begin{array}{ll}
	2&\hbox{if $2\nodiv d$ or $\mu_4\not\subset k$}\\
	3&\hbox{if $2|d$ and $\mu_4\subset k$}
\end{array} \right.
+\left\{ \begin{array}{ll}
	1&\hbox{if $3\nodiv d$}\\
	2&\hbox{if $3|d$ and $\mu_3\not\subset k$}\\
	3&\hbox{if $3|d$ and $\mu_3\subset k$.}
\end{array}\right.
$$
\enddemo

5.7.
We only used the trivial half of part 1 of the lemma.  We included the
other half because it is of use if one wants to analyze $\varphi_d$ using
the rational map $F_d\ratto W_d$.

\vglue-8pt
\section{The Tate conjecture}\label{s:Tate}
\vglue-4pt

From now on, we take $k$ to be $\Fq$, the field of $q$ elements, where
$q$ is a power of~$p$.

\demo{{\rm 6.1}} 
Let $X$ be a variety over $\Fq$.  The zeta function of $X$ is by
definition 
$$
\zeta(X,s)=\prod_{x}\frac{1}{1-q^{-\deg(x)s}}
$$
where the product is over all closed points of $X$ and $\deg(x)$ is the
degree of the residue field at $x$ as an extension of $\Fq$.  (For our
purposes it is enough to view this as a formal series in $q^{-s}$.)  
Alternatively, $\zeta(X,s)=Z(X,q^{-s})$ where
$$
Z(X,T)=\exp\left(\sum_{n=1}^\infty\frac{N_nT^n}{n}\right)
$$
and $N_n$ is the number of points on $X$ rational over $\F_{q^n}$.

It is immediate from the definition that the zeta function is
multiplicative for disjoint unions: if $X=X_1\cup X_2$ is a disjoint
union, then $\zeta(X,s)=\zeta(X_1,s)\zeta(X_2,s)$.

It is a theorem of Dwork that $Z(X,T)$ is a rational function of
$T$.  We will need the deeper connection with cohomology.  Fix a prime
$\ell\neq p$ and write $H^i(X)$ for the $\ell$-adic \'etale cohomology
group $H^i_{\hbox{\ninerm \'et}}(X\times_{\spec \Fq}\spec\Fqbar,\Q_\ell)$.  Then 
\begin{equation}\label{eqn:ZetaCohom}
Z(X,T)=\prod_{i=0}^{2\dim X}\det\left(1-Fr^*\,T|H^i(X)\right)^{(-1)^{i+1}} \speqnu{6.1.1}
\end{equation}
where $Fr:X\to X$ is the $q$-power Frobenius endomorphism.

Deligne's proof of the Weil conjectures implies that when $X$ is
smooth and proper over $\Fq$, the eigenvalues of $Fr^*$ on $H^i(X)$ are
algebraic numbers independent of $\ell$ and have absolute value
$q^{i/2}$ in any complex embedding.
\enddemo

6.2.
Now assume that $X$ is a smooth and proper variety over $\Fq$.  The
N\'eron-Severi group $NS(X)$ is defined to be the group of divisors on
$X$ modulo algebraic equivalence over $\Fqbar$.  
There is a cycle class map $NS(X)\to
H^2(X)$ which induces an injection $NS(X)\tensor\Q_\ell\to
H^2(X)^{Fr=q}$ where the exponent signifies the subspace where $Fr^*$
acts by multiplication by $q$.  This, together with the cohomological 
description of zeta functions, gives inequalities
\begin{equation}\label{eq:rk-ord}
\rk NS(X)\le \dim_{\Q_\ell}H^2(X)^{Fr=q}\le -\ord_{s=1}\zeta(X,s). \speqnu{6.2.1}
\end{equation}
 
The Tate conjecture (\cite{Tat65}) asserts that these are all
equalities.  We will refer to this assertion as ``(T) for $X$.''  (It
would be more precise to refer to this as (T1), since there are
conjectures for cycles of every codimension, but we will not need the
others.  Also, there are refined conjectures relating the leading
coefficient of the Taylor expansion of $\zeta(X,s)$ at $s=1$ to other
invariants of~$X$.)

\demo{{\rm 6.3}}
In the case where $X\to\Curve$ is an elliptic surface over $\Fq$ with
a section (so that the generic fiber is an elliptic curve $E$ over the
function field $K=\Fq(\Curve)$), the Tate conjecture for $X$ is
equivalent to the Birch and Swinnerton-Dyer conjecture for $E$.  (More
precisely (T) for $X$ implies that $\rk E(K)=\ord_{s=1}L(E/K,s)$.
Moreover, Tate proved (\cite{Tat66}) that when (T) holds the refined
conjecture of Birch and Swinnerton-Dyer on the leading Taylor
coefficient of $L(E/K,s)$ is true up to a power of $p$.  Milne showed
(\cite{Mil75}) that the full refined conjecture is true.  We only need the rank conjecture.)

Still assuming that $X$ is an elliptic surface, the cohomological
expression~\ref{eqn:ZetaCohom} for the zeta function and the Euler
characteristic formula of Groth\-endieck, Ogg, and Shafarevitch lead to
an upper bound on $\ord_{s=1}L(E/K,s)$ and thus also to an upper bound
on $\rk E(K)$.  If $g$ is the genus of $\Curve$ and $\n$ is the
conductor of $E$ (an effective divisor on $\Curve$), then we have
$$
\ord_{s=1}L(E/K,s)\le 4g-4+\deg(\n)
$$
if $E/K$ is nonconstant, and
$$
\ord_{s=1}L(E/K,s)\le 4g
$$
if $E/K$ is constant.  These bounds are ``geometric'' in that they are
insensitive to the finite field $\Fq$.  As we will explain in
Section~\ref{s:Bounds}, there exists a more refined arithmetic bound,
and asymptotically this bound is met by the curves in
Theorem~\ref{thm:intro}.

\specialnumber{6.4} \proclaim{Proposition}\label{prop:Tate}
Let $E$ be the elliptic curve over $\Fq(t)$ defined by
Equation~{\rm \ref{eq:E},} let $\EE$ be the elliptic surface over $\Fq$
defined in Section~{\rm \ref{s:EE},} and let $F_d/\Gamma$ be the quotient of
the Fermat surface of degree $d$ over $\Fq$ defined in
Section~{\rm \ref{s:Fd}.} 
\begin{itemize}
\ritem{1.} The Tate conjecture holds for $\EE${\rm .}  Equivalently{\rm ,}
the conjecture of Birch and Swinnerton\/{\rm -}\/Dyer holds for $E${\rm .}
\ritem{2.} $\rk E(\Fq(t))=-\ord_{s=1}\zeta(F_d/\Gamma,s)-1
+\varepsilon$ where 
$$
\varepsilon=\left\{ \begin{array}{ll}
	0&\hbox{if $2\nodiv d$ or $4\nodiv q-1$}\\
	1&\hbox{if $2|d$ and $4|q-1$}
\end{array}\right.
+\left\{ \begin{array}{ll}
	0&\hbox{if $3\nodiv d$}\\
	1&\hbox{if $3|d$ and $3\nodiv q-1$}\\
	2&\hbox{if $3|d$ and $3|q-1$.}
\end{array}\right.
$$
\end{itemize}
\endproclaim

  6.5. {\it Remark}.
The $\varepsilon$ term accounts for the points on
$E$ with either $x=0$ or $y=0$.

\demo{Proof}
Part 1 follows from the existence of a dominant rational map
$F_d\ratto \EE$ and well-known results on the Tate conjecture.  (We
refer to \cite{Tat94}, especially \S 5 for these results.)
Indeed, if $X\ratto Y$ is a dominant rational map, (T) for $X$ implies
(T) for $Y$.  But (T) is trivial for curves, and its truth for two
varieties implies it for their product.  Since a Fermat variety is
dominated by a product of curves (\cite{SK79}), (T) follows
for Fermat varieties, and this implies (T) for $\EE$.

The equivalence of (T) for $\EE$ and the conjecture of Birch and
Swinnerton-Dyer for $E$ was already noted above.  

To prove 2, we use the geometric analysis of Section~\ref{s:geometry}
and the multiplicativity of zeta functions.  From the Shioda-Tate
formula (Subsection~3.5) and Part 1 we have 
\begin{eqnarray*}
\rk E(K)&=&\rk NS(\EE)-\rk L\\
&=&-\ord_{s=1}\zeta(\EE,s)-\rk L. 
\end{eqnarray*}
Since $\EE$ is obtained from $W$ by blowing up, i.e., by removing
points and adding curves, multiplicativity of zeta functions yields
that
\begin{eqnarray*}
-\ord_{s=1}\zeta(\EE,s)&=&-\ord_{s=1}\zeta(W,s)+\sum_v(n_v-1)\\
&=&-\ord_{s=1}\zeta(W,s)+\rk L -2
\end{eqnarray*}
where $n_v$ is the number of irreducible components in the fibre at
$v$.  Combining multiplicativity with the geometric analysis of
Section~\ref{s:geometry} we see that
$$
-\ord_{s=1}\zeta(W,s)-4-\varepsilon
=-\ord_{s=1}\zeta(F_d/\Gamma,s)-3.
$$
Assembling these
ingredients gives the desired formula.\enddemo

6.6.
Note that we did not need the actual values of $n_v$ in the proof,
since they cancel out.  Nevertheless, we recorded them in
Section~\ref{s:E} for future use, e.g., for height computations.
\pagebreak

\section{The zeta function of a Fermat surface}

7.1.
The zeta functions of Fermat varieties were computed in terms of
Gauss and Jacobi sums by Weil in his landmark paper \cite{Wei49}.
We will need a refinement of this calculation due to Shioda which
takes into account the action of $G$, and we will need to make the
relevant Jacobi sums explicit.  Remarkably, an explicit calculation of
zeta functions of Fermat varieties \emph{over the prime field} does
not seem to be in the literature.  When doing explicit calculations,
most authors pass immediately to the case where $q\equiv1\hbox{ mod}\, d$.  In
the next two sections, we will explicitly compute the zeta function of
$F_d/\Gamma$ over any finite field $\Fq$, in the ``supersingular''
case, i.e., when $d$ divides $p^n+1$ for some positive integer $n$.

\demo{{\rm 7.2}} 
We will use the cohomological description
$$
\zeta(X,s)=Z(X,q^{-s})=\prod_{i=0}^{2\dim X}P_i(X,q^{-s})^{(-1)^{i+1}}
$$
where
$$
P_i(X,T)=\det\left(1-Fr^*\,T|H^i(X)\right).
$$

Since $H^i(F_d/\Gamma)=H^i(F_d)^{\Gamma}$ and the eigenvalues of
Frobenius on $H^i(F_d)$ have absolute value $q^{i/2}$, only
$P_2(F_d/\Gamma,q^{-s})$ can contribute to the order of pole of
$\zeta(F_d/\Gamma,s)$ at $s=1$.  Thus we will concentrate on
$H^2(F_d)$ and its $\Gamma$ invariants.
\enddemo

7.3.
Fix an algebraic closure $\Qbar$ of $\Q$.  Let
$\p$ be a prime of $\O_\Qbar$, the ring of integers of $\Qbar$, over
$p$.  We view all finite fields of characteristic $p$ as subfields of
$\O_\Qbar/\p$, which is an algebraic closure of $\Fp$.

Reduction modulo $\p$ induces an isomorphism between the group of all
roots of unity of order prime to $p$ in $\O_\Qbar$ and the
multiplicative group of $\O_\Qbar/\p$.  We let
$t:(\O_\Qbar/\p)^\times\to\Qbar^\times$ denote the inverse of this
isomorphism.  We will use the same letter $t$ for the restriction to
any finite field $\Fq^\times$.

Fix an algebraic closure $\Qbar_\ell$ of $\Q_\ell$ and an embedding
$\Qbar\into\Qbar_\ell$.  For convenience, we will assume that
$\ell\equiv1\, {\rm mod}\, {pd}$ so that $\Q_\ell$ contains all the $pd^{\rm th}$ roots
of unity.

\demo{{\rm 7.4}} 
Now we introduce Gauss and Jacobi sums.  Fix a nontrivial character
$\psi_0:\Fp\to\Qbar^\times$ and for each finite extension $\Fpf$ of
$\Fp$, let $\psi:\Fpf\to\Qbar^\times$ be defined by
$\psi=\psi_0\compose\tr_{\Fpf/\Fp}$.  If
$\chi:\Fpf^\times\to\Qbar^\times$ is a nontrivial character, we
define a Gauss sum by
$$
g(\chi,\psi)=-\sum_{x\in\Fpf^\times}\chi(x)\psi(x).
$$

If $\chi_1,\dots,\chi_n$ are characters
$\Fpf^\times\to\Qbar^\times$, not all trivial, such that the product
$\chi_1\cdots\chi_n$ is trivial, we define a Jacobi sum by
$$
J(\chi_1,\dots,\chi_n)=\frac1{p^f-1}
\sum_{{x_1,\dots,x_n\in\Fpf^\times\atop x_1+\cdots+x_n=0}}
\chi_1(x_1)\cdots\chi_n(x_n).
$$

It is well-known (see \cite{Wei49} for example) that 
$$
J(\chi_1,\dots,\chi_n)=\left\{ \begin{array}{ll}
\frac{(-1)^n}{p^f}\prod_{i=1}^n g(\chi_i,\psi)&\hbox{if all $\chi_i$ are nontrivial}\\
0&\hbox{otherwise}
\end{array}\right.
$$
and that (in any complex embedding)
$|g(\chi,\psi)|=p^{f/2}$.  
\enddemo

7.5.
Recall the group $G$ of automorphisms of $F_d$ introduced in
Section~\ref{s:Fd}.  Let $\hat G$ denote the group of characters $G$
with values in $\Qbar$ (and thus also $\Qbar_\ell$ via our fixed
embedding).  Using the character
$t:(\O_\Qbar/\p)^\times\to\Qbar^\times$, we can identify $\hat G$ with
$$
\left\{a=(a_0,a_1,a_2,a_3)\in (\Z/d\Z)^4\left|\sum a_i=0\right.\right\}
$$
where the pairing
$G\times\hat G\to\Qbar^\times$ is
$$
a(z)=
\left\langle(a_0,a_1,a_2,a_3),[\zeta_0,\zeta_1,\zeta_2,\zeta_3]\right\rangle
=\prod_{i=0}^3t(\zeta_i)^{a_i}.
$$
Note that by our assumption that $\ell\equiv1\, {\rm mod}\, {pd}$, the values of
$a\in\hat G$ lie in $\Q(\mu_d)\subset\Q_\ell$.

For $a\in \hat G$, we denote by $H^2(F_d)(a)$ the subspace of classes
$c\in H^2(F_d)$ such that $z^*(c)=a(z)c$ for all $z\in G$.  Also, we
write $qa$ for $(qa_0,\dots,qa_3)$ if $a=(a_0,\dots,a_3)$.

Recall that $Fr:F_d\to F_d$ denotes the $q$-power Frobenius
endomorphism.  From the formula
$$
Fr\compose[\zeta_0,\zeta_1,\zeta_2,\zeta_3]=
[\zeta_0^q,\zeta_1^q,\zeta_2^q,\zeta_3^q]\compose Fr
$$
we deduce that $Fr^*$ sends $H^2(F_d)(a)$ to $H^2(F_d)(qa)$.  For each
$a\in \hat G$, we let $u(a)$ denote the smallest positive integer such
that $q^{u(a)}a=a$.  Then $(Fr^{u(a)})^*$ maps $H^2(F_d)(a)$ to
itself.

It turns out that $H^2(F_d)(a)$ is zero or 1-dimensional, and
$(Fr^{u(a)})^*$ acts by multiplication by a Jacobi sum.  More precisely,
for $a\in \hat G$, $a\neq0$, define a
Jacobi sum $J(a)$ as follows: Let
$\chi_i:\F_{q^{u(a)}}^\times\to\Qbar$ be defined as
$\chi_i=t^{\frac{q^{u(a)}-1}{d}a_i}$ and set
$J(a)=J(\chi_0,\dots,\chi_3)$.  Note that $J(qa)=J(a)$.  By
convention, we set $J(0)=q$.\pagebreak

\proclaimtitle{Shioda}
\specialnumber{7.6} \proclaim{Proposition} 
Let $F_d$ be the Fermat surface of degree $d$ over $\Fq$ and let $\hat
G'=\{a=(a_0,\dots,a_3)\in\hat G\ |\ a=0\hbox{ or }a_i\neq0\hbox{
for }i=0,\dots3\}${\rm .}
\smallbreak
{\rm 1.} $H^2(F_d)(a)$ is zero if $a\not\in\hat
G'$ and is $1$\/{\rm -}\/dimensional if $a\in\hat G'${\rm .}
\smallbreak
{\rm 2.} If $a\in\hat G'$ then $(Fr^{u(a)})^*$ acts on
$H^2(F_d)(a)$ by multiplication by $J(a)${\rm .  }
\smallbreak
 {\rm 3.} \hangindent=36pt\hangafter=1 If $a\in\hat G'$ then the characteristic polynomial of $Fr^*$ on
$\bigoplus_{i=0}^{u(a)-1}H^2(F_d)(q^ia)$ is equal to
$(1-J(a)T^{u(a)})${\rm . }

\endproclaim

\demo{Proof}
In~\cite{SK79}, Shioda and Katsura show that the cohomology
of a Fermat variety is built up from the cohomology of lower
dimensional Fermat varieties of the same degree.  This allows one to
reduce to the case of curves.  The proposition for Fermat curves
is~\cite[Cor.~2.4]{Kat81}.

We will sketch another proof, closely related to that in
\cite{Kat81}, which works uniformly in all dimensions.
For simplicity we will only discuss the case of Fermat surfaces.

To that end, consider the finite morphism $\pi:F_d\to
F_d/G\cong\P^2$.  The sheaf $\FF=\pi_*\Q_\ell$ carries a natural
action of $G$ and we have a decomposition $\FF=\mathbold{\oplus}_{a\in\hat
G}\FF(a)$.  By the Leray spectral sequence for $\pi$,
$$H^2(F_d)(a)=H^2_{\hbox{\ninerm \'et}} \left(\P^2\times\spec\Fqbar,\FF\right)(a)
=H^2_{\hbox{\ninerm \'et}}\left(\P^2\times\spec\Fqbar,\FF(a)\right).$$ 

Now each $\FF(a)$ is lisse of rank 1 on the locus
$\{[y_0,\dots,y_3]|y_i\neq0$ if\break $a_i\neq0\}\subset\P^2$ and is
zero elsewhere.  If $r=q^f$ is a power of $q^{u(a)}$ and $y$ is an
$\F_{r}$-rational point of $\P^2$ with corresponding Frobenius
$Fr_y=Fr^f$, then $Fr_y$ induces an automorphism of the geometric
stalk, $Fr_y:\FF(a)_{\overline{y}}\to\FF(a)_{\overline{y}}$.  We have  
$\tr Fr_y|\FF(a)_{\overline{y}}=\chi_{a,r}(y)$ where $\chi_{a,r}$ is
defined by
\begin{eqnarray*}
\chi_{a,r}([y_0,\dots,y_3])
&=&\prod_{i=0}^3t^{\frac{r-1}da_i}(y_i)\\ 
&=&\prod_{i=0}^3t^{\frac{q^{u(a)}-1}da_i}
\left(N_{\F_r/\F_{q^{u(a)}}}(y_i)\right)
\end{eqnarray*}
and we interpret $t^b(0)$ as 0 if $b\neq0$ and as 1 if $b=0$.

Using the Grothendieck-Lefschetz trace formula and the Hasse-Davenport
relation, we see that
\begin{eqnarray*}
&&\prod_{i=0}^4\det\left(1-(Fr^{u(a)})^*T|H^i(F_d)(a)\right)^{(-1)^{i+1}}\\
&&\qquad =
\left\{ \begin{array}{ll}
(1-J(a)T)^{-1}&\hbox{if $a\neq0$}\\ 
(1-T)^{-1}(1-qT)^{-1}(1-q^2T)^{-1}&\hbox{if $a=0$}.
\end{array}\right.
\end{eqnarray*}
\pagebreak 

\noindent
 Now Deligne's purity theorem for the $H^i(F_d)$ implies that
$H^2(F_d)(a)$ is either 0 or 1-dimensional, and is nonzero if and
only if $J(a)\neq0$, i.e., if and only if $a\in\hat G'$.  We also see
that $(Fr^{u(a)})^*$ acts on $H^2(F_d)(a)$ by multiplication by~$J(a)$.

Part~3 is an easy consequence of Part~2.\enddemo

The proposition reduces the problem of computing the order of pole of
the zeta function of $F_d/\Gamma$ to computing some Jacobi sums.  
Let $\Gamma^\perp\subset\hat G$ be the set of
characters which are trivial on $\Gamma\subset G$.  Clearly
$\Gamma^\perp$ is the cyclic subgroup of order $d$ generated by
$(3,-6,2,1)$.  With this notation, we have:

\specialnumber{7.7} \proclaim{{C}orollary}\label{cor:ord-orbits}
Let $A_1,\dots,A_k$ be the orbits of multiplication by $q$ on
$\Gamma^\perp\cap\hat G'$ and choose $a_i\in A_i${\rm .}  Then
$$
P_2(F_d/\Gamma,T)=\prod_{i=1}^k(1-J(a_i)T^{u(a_i)}).
$$
\endproclaim
\vglue-9pt
 
\section{Explicit Gauss and Jacobi sums}
\vglue-6pt
\specialnumber{8.1} \proclaim{Proposition}\label{prop:ExJac}
Suppose that some power of $p$ is congruent to $-1$ modulo
$d${\rm .}  Then for all $a\in\hat G'${\rm ,} $J(a)=q^{u(a)}${\rm .}
\endproclaim

The rest of this section is devoted to proving the proposition.

\specialnumber{8.2} \proclaim{Lemma}\label{lemma:denoms}
Let $b$ be a rational number with $0<b<1$ and suppose that there exist
positive integers $n$ and $f$ such that $(p^n+1)b\in\Z$ and
$(p^f-1)b\in\Z${\rm .}  If $b\neq1/2$ then $f$ is even and setting
$e=\gcd(n,f/2)${\rm ,} we have $(p^e+1)b\in\Z${\rm .}
\endproclaim

\demo{Proof}
Write $b=c/d$ in lowest terms.  If $b\neq1/2$ then $d>2$.  Our
assumptions imply that $p^n\equiv-1\, {\rm mod}\,  d$ and $p^f\equiv1\, {\rm mod}\,  d$, and
so $f$ must be even since $d>2$.  Now
$$
\gcd(p^{2n}-1,p^f-1)=p^{\gcd(2n,f)}-1=p^{2e}-1
$$
so $p^{2e}\equiv1\, {\rm mod}\,  d$.  Since $p^n=(p^e)^{(n/e)}\equiv-1\, {\rm mod}\,  d$, we
must have that $n/e$ is odd and $p^e\equiv-1\, {\rm mod}\,  d$.  Thus
$(p^e+1)b\in\Z$, as desired.\enddemo

\proclaimtitle{Shafarevitch-Tate}
\specialnumber{8.3} \proclaim{Lemma}\label{lemma:ST}
 Let $\chi:\F_{p^{2f}}^\times\to\Qbar^\times$ be a
nontrivial character which is trivial on $\F_{p^f}^\times${\rm .}  Then
$g(\chi,\psi)=-\chi(x)p^f$ where $x\in\F_{p^{2f}}^\times$ is any
element with $\tr_{\F_{p^{2f}}/\F_{p^{f}}}(x)=0${\rm .}
\endproclaim

\demo{Proof}
Recall that $\psi(x)=\psi_0(\tr_{\F_{p^{2f}}/\Fp}(x))$ where $\psi_0$
is a fixed nontrivial character of $\Fp$.  Abbreviating
$\tr_{\F_{p^{2f}}/\F_{p^{f}}}$ to $\tr$, we have \pagebreak

\begin{eqnarray*}
\noalign{\vskip-15pt}
g(\chi,\psi)&=&-\sum_{x\in\F_{p^{2f}}^\times}\chi(x)\psi(x)\\[5pt]
&=&-\sum_{x\in\F_{p^{2f}}^\times/\F_{p^f}^\times}\chi(x)
    \sum_{y\in\F_{p^{f}}^\times}\psi(xy)\\[5pt]
&=&-\sum_{x\in\F_{p^{2f}}^\times/\F_{p^f}^\times}\chi(x)
    \left\{ \begin{array}{rl}
    p^f-1& \hbox{if $\tr(x)=0$}\\[5pt]
    -1 & \hbox{if $\tr(x)\neq0$}
    \end{array}\right.\\[5pt]
&=&-\sum_{{x\in\F_{p^{2f}}^\times/\F_{p^f}^\times\atop\tr(x)=0}}
\chi(x)p^f.
\end{eqnarray*}
But $\tr:\F_{p^{2f}}\to\F_{p^f}$ is $\F_{p^f}$-linear and surjective,
so its kernel is a 1-dimensional $\F_{p^f}$-vector space.  This means
that there is just one term in the last displayed sum, and this proves
the lemma.\enddemo

\demo{Proof of Proposition~{\rm \ref{prop:ExJac}}}
If $a=0$ then $u(a)=1$ and $J(a)=q$ by definition.  
If $a=(d/2,d/2,d/2,d/2)$, then $u(a)=1$ and
$J(a)=g(t^{(q-1)/2},\psi)^4/q$.  But it is elementary and well-known
that $g(t^{(q-1)/2},\psi)=\pm\sqrt{\pm q}$, and so $J(a)=q$.

Now assume that $a\in\hat G'$ and $a\neq0$, $a\neq(d/2,d/2,d/2,d/2)$.
Lemma~\ref{lemma:denoms}, applied to the $b_i=a_i/d$, shows that $u(a)$
is even.  Setting $g_i=g(t^{\frac{q^{u(a)}-1}da_i},\psi)$ and
$e=\gcd(n,u(a)/2)$ we have
\begin{eqnarray*}
g_i&=&-\sum_{x\in\F_{q^{u(a)}}^\times}t^{\frac{q^{u(a)}-1}da_i}(x)\psi(x)\\
&=&-\sum_{x\in\F_{q^{u(a)}}^\times}t^{\frac{q^{2e}-1}da_i}
\left({\rm N}_{\F_{q^{u(a)}}/\F_{q^{2e}}}(x)\right)\psi(x)\\
\noalign{\noindent which, by the Hasse-Davenport relation, is}
&=&\left(-\sum_{x\in\F_{q^{2e}}^\times}t^{\frac{q^{2e}-1}da_i}(x)
\psi(x)\right)^{u(a)/(2e)}.
\end{eqnarray*}
(We have abusively written $\psi$ for the additive
characters of both $\F_{q^{u(a)}}$ and $\F_{q^{2e}}$.)  But by
Lemma~\ref{lemma:denoms}, $(q^e+1)a_i/d\in\Z$ and so the inner sum is
of   \pagebreak the type considered in Lemma~\ref{lemma:ST}.  Thus
$$
g_i=\left(-t^{\frac{q^{2e}-1}da_i}(x)\right)^{u(a)/(2e)}q^{u(a)/2}
$$
where $x\in\F_{q^{2e}}^\times$ is any element with
$\tr_{\F_{q^{2e}}/\F_{q^{e}}}(x)=0$. 

Taking the product over $i=0,\dots,3$ and using the fact that
$\sum_{i=0}^3a_i\equiv0\, {\rm mod}\,  d$, we see that $J(a)=q^{u(a)}$.\enddemo
 
\section{The rank of $E$}

9.1.
We are now in  a position to compute the rank of $E_d$ over $\Fq(t)$
for any $d$ dividing $p^n+1$.
By Proposition~\ref{prop:Tate}, the rank is
$$
-\ord_{s=1}\zeta(F_d/\Gamma,s)-1+
\left\{ \begin{array}{ll}
	0&\hbox{if $2\nodiv d$ or $4\nodiv q-1$}\\
	1&\hbox{if $2|d$ and $4|q-1$}
\end{array}\right.
+\left\{ \begin{array}{ll}
	0&\hbox{if $3\nodiv d$}\\
	1&\hbox{if $3|d$ and $3\nodiv q-1$}\\
	2&\hbox{if $3|d$ and $3|q-1$.}
\end{array}\right.
$$

By Corollary~\ref{cor:ord-orbits} and Proposition~\ref{prop:ExJac},
$-\ord_{s=1}\zeta(F_d/\Gamma,s)$ is equal to the number of orbits of
multiplication by $q$ on $\Gamma^\perp\cap\hat G'$.  Here 
$$
\hat G'=\{a\in\hat G\ |\
a=0\hbox{ or }a=(a_0,\dots,a_3)\hbox{ with }a_i\neq0\}
$$
and
$\Gamma^\perp$ is the cyclic subgroup of $\hat G$ generated by
$(3,-6,2,1)$.  Thus $\Gamma^\perp\cap\hat G'$ is in bijection with
$\{a_3\in\Z/d\Z\ |\ 6a_3\neq0\}\cup\{0\}$.  The size of the orbit of
$a_3$ depends only on $e=\frac d{(d,a_3)}$ and is equal to $o_e(q)$, the order
of $q$ in $(\Z/e\Z)^\times$.  Thus we have
$$
-\ord_{s=1}\zeta(F_d/\Gamma,s)=
\sum_{{e|d\atop e\nodiv6}}\frac{\phi(e)}{o_e(q)}+1.
$$
(The term 1 corresponds to $e=1$, i.e., to $a_3=0$.)

Putting everything together, we have our main theorem.

\specialnumber{9.2} \proclaim{Theorem}\label{thm:main}
Let $p$ be a prime{\rm ,} $n$ a positive integer{\rm ,} and $d$ a divisor of
$p^n+1${\rm .}  Let $q$ be a power of $p$ and let $E$ be the elliptic curve
over $\Fq(t)$ defined by
$$
y^2+xy=x^3-t^d.
$$
Then the $j$\/{\rm -}\/invariant of $E$ is not in $\Fq${\rm ,} the conjecture of Birch
and Swinnerton\/{\rm -}\/Dyer holds for $E${\rm ,} and the rank of $E(\Fq(t))$ is
$$
\sum_{ e|d\atop e\nodiv6}\frac{\phi(e)}{o_e(q)}
+\left\{ \begin{array}{ll}
	0&\hbox{if $2\nodiv d$ or $4\nodiv q-1$}\\
	1&\hbox{if $2|d$ and $4|q-1$}
\end{array}\right.
+\left\{ \begin{array}{ll}
	0&\hbox{if $3\nodiv d$}\\
	1&\hbox{if $3|d$ and $3\nodiv q-1$}\\
	2&\hbox{if $3|d$ and $3|q-1$.}
\end{array}\right.
$$
Here $\phi(e)$ is the cardinality of $(\Z/e\Z)^\times$ and $o_e(q)$
is the order of $q$ in $(\Z/e\Z)^\times$.
\endproclaim

9.3.
We now specialize to the case where $d=p^n+1$.  If $q=p$ then
$o_e(p)\le 2n$ for all divisors $e$ of $d$.  Applying the theorem, we
see that the rank of $E$ over $\Fp(t)$ is at least $(p^n-1)/2n$.
This completes the proof of Theorem~\ref{thm:intro}.

On the other hand, if we take $q$ to be a power of $p^{2n}$, then
$o_e(q)=1$ for all divisors $e$ of $d$. The theorem then implies that
the rank of $E$ over $\Fq(t)$ is $d-1=p^n$ if $6\nodiv d$ and
$d-3=p^n-2$ if $6|d$.

\section{Rank bounds}\label{s:Bounds}

10.1.
Let $\Curve$ be a smooth complete curve of genus $g$ over $\Fq$ and
let $E$ be an elliptic curve over $K=\Fq(\Curve)$.  Write $\n$ for the
conductor of $E$ and let $\deg(\n)$ be the degree of $\n$, viewed as
an effective divisor on $\Curve$.
As mentioned in Section~\ref{s:Tate}, there is a
geometric bound on the rank of $E$:
$$
\rk E(K)\le\ord_{s=1}L(E/K,s)\le 
\left\{ \begin{array}{ll}4g&\hbox{if $E$ is constant}\\
4g-4+\deg(\n)&\hbox{if $E$ is not constant.}\end{array}\right.
$$
This bound is geometric in that it is not affected if we extend the
constant field $\Fq$.  But as we have seen in the previous section,
both $\rk E(K)$ and $\ord_{s=1}L(E/K,s)$ can change dramatically if the
constant field is enlarged.

In fact, there is an arithmetic bound, i.e., a bound which is
sensitive to the finite field of constants.  In
\cite[Prop.~6.9]{Bru92} Brumer used Weil's ``explicit formula''
technique to prove an upper bound
$$
\ord_{s=1}L(E/K,s)\le\frac{4g-4+\deg(\n)}{2\log_q\deg(\n)}
+C\frac{\deg(\n)}{(\log_q\deg(\n))^2}
$$
where $\log_q$ denotes the logarithm to base $q$ and $C$ is an
explicit constant depending only on $g$ and $q$.  (Here we ignore the
finitely many elliptic curves over $K$ with trivial conductor.)  This
is the function field analogue of a theorem of Mestre \cite{Mes86}
which says that if $E$ is a modular elliptic curve over $\Q$ then,
assuming a generalized Riemann hypothesis, $$\ord_{s=1}L(E/\Q,s)=O(\log
N/\log\log N).$$

Brumer's bound is visibly sensitive to the field of constants and is an
improvement on the geometric bound when $\deg(\n)$ is
large with respect to $q$.  

\demo{{\rm 10.2}}
The curves of Theorem~\ref{thm:main} show that the main term of
Brumer's arithmetic bound over $\Fp(t)$ is sharp.  Indeed, the curve
with $d=p^n+1$ has $\deg(\n)=p^n+4$ if $6\nodiv d$ and
$\deg(\n)=p^n+2$ if $6|d$, and its rank (analytic and algebraic) is at
least $(p^n-1)/2n$.  (Note also that these curves meet the geometric
bound over $\Fq(t)$ when $\Fq$ contains $\F_{p^{2n}}$.)
\enddemo

10.3.
The (isotrivial) elliptic curves of \cite{TS67} also
meet the main term of Brumer's arithmetic bound over $\Fp(t)$.
Indeed, the curve of their Theorem~2 (with $f=p^n+1$) has $\deg(\n)$
approximately $2p^n$ and rank approximately $p^n/n$.

\demo{{\rm 10.4}}
The preceding remarks show that the Brumer arithmetic bound is
asymptotically sharp in the function field case.  We believe that the
Mestre bound should likewise be asymptotically sharp in the number
field case.

Let $K$ now be a number field.  For a positive integer $N$, let
$r_K(N)$ be the maximum, over all elliptic curves $E$ over $K$ with
conductor $\n$ satisfying $N_{K/\Q}(\n)=N$, of $\rk E(K)$; if there
are no such curves, we set $r_K(N)=0$.  Assuming various standard
conjectures, it follows from a simple generalization of Mestre's
argument that $r_K(N)=O(\log N/\log\log N)$ (where the constant of
course depends on $K$), and so the limit in the following conjecture
is finite.
\enddemo

\specialnumber{10.5}\proclaim{{C}onjecture}
$$
\limsup_N\frac{r_K(N)}{\log N/\log\log N}>0.
$$
\endproclaim

10.6.
If $E$ is an elliptic curve over $\Q$, let $N_\Q(E)$ be its conductor
and let $N_K(E)$ be the norm from $K$ to $\Q$ of the conductor of $E$
viewed as elliptic curve over $K$.  Then there is a constant
$C$ depending only on $K$ such that
$$
1\le\frac{N_\Q(E)^{[K:\Q]}}{N_K(E)}\le C
$$
for all elliptic curves $E$ over $\Q$.
This can be used to prove that the conjecture
for a general number field $K$ follows from the conjecture for $\Q$.
\vglue12pt
10.7 (Added in proof).  If one replaces ``elliptic curve" with ``abelian variety" in the definition of $r_K(N)$ then the limit
in Conjecture 10.5 is still finite and the conjecture has been proven unconditionally.  See E. Kowalski and\break P. Michel,
{\it Acta Arith.\/} {\bf 94} (2000), 303--343.

\AuthorRefNames [MM00]


\begin{references}
 \bibitem{Bru92} 
\name{A.\ Brumer}, The average rank of elliptic curves. I,
 {\it Invent.\ Math\/}.
  {\bf 109} (1992), 445--472. 


\bibitem{GZ86} \name{B.~H. Gross} and \name{D.~B. Zagier}, Heegner points and derivatives
of
  ${L}$-series, {\it Invent.\ Math\/}.\  {\bf 84} (1986), 225--320.
 

\bibitem{Kat81} 
\name{N.~M. Katz}, Crystalline cohomology, {D}ieudonn\'e modules, and
{J}acobi
  sums, in {\it Automorphic Forms{\rm ,} Representation Theory and Arithmetic\/}
(Bombay,
  1979), Tata Inst. Fundamental Res., Bombay, 1981, 165--246.
 

\bibitem{Mes86} 
\name{J.-F. Mestre}, Formules explicites et minorations de conducteurs de
  vari\'et\'es\break alg\'ebriques, {\it Compositio Math\/}.\  {\bf 58}
(1986), 209--232. 

\bibitem{Mil75} 
\name{J.~S. Milne}, On a conjecture of {A}rtin and {T}ate, {\it Ann.\ of
Math\/}.\
  {\bf 102} (1975), 517--533.  See also the addendum available at http://www.jmilne.org/math/ 

\bibitem{MM00} 
\name{R.~Martin} and \name{W.~McMillen}, An elliptic curve over ${\bf Q}$ with rank at
least
  24, Preprint (2000), Posted to the Usenet newsgroup math.sci.nmbrthry by
  V.~Miller on May 2, 2000.\ Available at  
  http://listserv.nodak.edu/archives/nmbrthry.html. 

\bibitem{MP86} 
\name{I.~Morrison} and \name{U.~Persson}, Numerical sections on elliptic
surfaces,
  {\it Compositio Math\/}.\ {\bf 59} (1986), 323--337. 

\bibitem{Shi72} 
\name{T.~Shioda}, On elliptic modular surfaces, {\it J.\ Math.\ Soc.\ Japan\/}
{\bf 24} (1972), 20--59. 

\bibitem{Shi86} 
\bibline, An explicit algorithm for computing the {P}icard number
of certain algebraic surfaces, {\it Amer.\ J.\ Math\/}.\ {\bf 108}
(1986), 415--432. 

\bibitem{SK79} 
\name{T.~Shioda} and \name{T.~Katsura}, On {F}ermat varieties, {\it T{\rm \^{\it o}}hoku Math.\
J\/}.\ {\bf 31} (1979), 97--115. 

\bibitem{Tat65} 
\name{J.~T. Tate}, Algebraic cycles and poles of zeta functions,
in {\it Arithmetical  Algebraic Geometry\/}, 93--110 (Proc.\ Conf.\ Purdue
Univ., 1963), Harper
\& Row, New York,
  1965. 

\bibitem{Tat66} 
\bibline, On the conjectures of {B}irch and {S}winnerton-{D}yer
and a
  geometric analog, {\it S\'eminaire Bourbaki\/}, {\bf 9}, Soc. Math.
France, Paris,
  1966, Exp.\ No.\ 306, 415--440, 1995. 

\bibitem{Tat75} 
\name{J. T. Tate}, Algorithm for determining the type of a singular fiber
in an
  elliptic pencil, in {\it Modular Functions of One Variable, IV\/}
(Proc.\ Internat.\
  Summer School, Univ. Antwerp, Antwerp, 1972), Springer-Verlag, New
York,  {\it Lecture Notes in Math\/}.\ {\bf 476} (1975), 33--52.

\bibitem{Tat94} 
\bibline, Conjectures on algebraic cycles in $l$-adic cohomology,
 in {\it  Motives\/} (Seattle, WA, 1991), {\it Proc.\ Sympos.\ Pure Math\/}.\ {\bf 55}, Amer. Math. Soc., Providence,
RI, 1994, 71--83. 

\bibitem{TS67} 
\name{J.~T. Tate} and \name{I.~R. {S}hafarevitch}, The rank of elliptic curves,
{\it   Akad.\ Nauk SSSR\/} {\bf175} (1967), 770--773. 

\bibitem{Wei49} 
\name{A.~Weil}, Numbers of solutions of equations in finite fields, {\it Bull.\
Amer.\ Math.\ Soc\/}.\ {\bf 55} (1949), 497--508.
\end{references}
\end{document}